# A Generalization of the Idea of Disjunction

Kerry M. Soileau

July 13, 2010

ABSTRACT

We generalize the concept of disjunction.

## 1. INTRODUCTION

Let $X$ be a nonempty set.

Definition 1: We say that a relation $\Diamond$ on $X$ is disjunctive if and only if for every $x_1, x_2 \in X$,

$$\forall y \in X \left( x_1 \Diamond y \leftrightarrow x_2 \Diamond y \right) \text{ if and only if } x_1 = x_2 \tag{1}$$

$$x_1 \Diamond x_2 \rightarrow x_2 \Diamond x_1 \tag{2}$$

$$\neg\left( x_1 \Diamond x_1 \right) \tag{3}$$

Examples: (1) Let $X$ be some nonempty set, and define $x_1 \Diamond x_2$ if and only if $x_1 \neq x_2$. (2) Take $X$ to be the power set of some nonempty set, and define $x_1 \Diamond x_2$ if and only if $x_1 \cap x_2 = \varnothing$. (3) Take $X = \{a, b, c\}$, and define $x_1 \Diamond x_2$ if and only if $\{x_1, x_2\} = \{a, b\}$.

Proposition 1: For any $x_1, x_2 \in X$, $x_1 \Diamond x_2 \rightarrow x_1 \neq x_2$.

Proof: Suppose for a contradiction that $x_1 \Diamond x_2$ and $x_1 = x_2$. This means $x_1 \Diamond x_1$, which contradicts (3). ■

## 2. A DISJUNCTIVE RELATION INDUCES A PARTIAL ORDER

Definition 2: Define the order $\prec$ satisfying $x_1 \prec x_2$ if and only if $\forall y \in X\left(x_2 \Diamond y \rightarrow x_1 \Diamond y\right)$.

Proposition 2: If $x_1 \Diamond x_2$, then $x_3 \prec x_2$ implies $x_1 \Diamond x_3$.

Proof: Suppose $x_1 \Diamond x_2$ and $x_3 \prec x_2$. Since $x_1 \Diamond x_2$, from (2) we get $x_2 \Diamond x_1$. Since $x_3 \prec x_2$, we have $\forall y \in X\left(x_2 \Diamond y \rightarrow x_3 \Diamond y\right)$ so in particular $x_2 \Diamond x_1 \rightarrow x_3 \Diamond x_1$ and thus $x_3 \Diamond x_1$. Another application of (2) yields $x_1 \Diamond x_3$. ■

Proposition 3: $\prec$ is a partial order on $X$.

Proof: (1) Fix $x \in X$. Since clearly $\forall y \in X\left(x \Diamond y \rightarrow x \Diamond y\right)$, it follows that $x \prec x$. (2) Suppose $x_1 \prec x_2$ and $x_2 \prec x_1$. Then $\forall y \in X\left(x_2 \Diamond y \rightarrow x_1 \Diamond y\right)$ and $\forall y \in X\left(x_1 \Diamond y \rightarrow x_2 \Diamond y\right)$, hence $\forall y \in X\left(x_1 \Diamond y \leftrightarrow x_2 \Diamond y\right)$, which by definition of $\Diamond$ implies $x_1 = x_2$. (3) Suppose $x_1 \prec x_2$ and $x_2 \prec x_3$. Fix $y \in X$. Then $x_2 \Diamond y \rightarrow x_1 \Diamond y$ and $x_3 \Diamond y \rightarrow x_2 \Diamond y$, implying $x_3 \Diamond y \rightarrow x_1 \Diamond y$, thus $x_1 \prec x_3$.■

Example:

Take $X = \{a,b,c\}$, and define $x_1 \Diamond x_2$ if and only if $\left(x_1, x_2\right) = (a,b)$ or $\left(x_1, x_2\right) = (b,a)$.

Then $x_1 \prec x_2$ if and only if $x_2 = c$.

This example can be illustrated with a table:

| ◊ | a | b | c |
|---|---|---|---|
| a | F | T | F |
| b | T | F | F |
| c | F | F | F |

The entry for each row and column gives the truth value of the disjunctive relation, for example the table indicates that $a \Diamond b$ is true and $b \Diamond c$ is false. When a table is constructed for a given disjunctive relation, it will be noticed that the table has three properties:

1. The diagonal contains only entries of $F$,
2. the table is symmetric with respect to transposing rows and columns, and
3. no two columns or are alike (and the same is true for rows).

International Space Station Program Office, Avionics and Software Office, Mail Code OD, NASA Johnson Space Center, Houston, TX 77058

E-mail address: ksoileau@yahoo.com